\documentclass{roalart}

\pdfoutput=1

\author[Roberto Albesiano]{Roberto Albesiano}
\address{Department of Pure Mathematics, University of Waterloo}
\email{roberto.albesiano@uwaterloo.ca}
\urladdr{https://math.uwaterloo.ca/~ralbesiano/}

\title{From division to extension}

\date{\today}

\begin{document}

\begin{abstract}
    We present a short proof of a version of the Ohsawa--Takegoshi--Manivel $L^2$ extension theorem as a corollary of a Skoda-type $L^2$ division theorem with bounded generators.  The new division theorem is of independent interest: the boundedness of generators allows to send the parameter $\alpha>1$ of the usual $L^2$ division theorems to 1 in the norm of the datum of the division.  As an aside, we also use the new division theorem to prove a Briançon--Skoda-type result.
\end{abstract}

\section*{Introduction}
The $L^2$ extension theorem and the $L^2$ division theorem are intimately related results in complex analytic geometry.  In fact, T.~Ohsawa proves in \cite{Ohsawa2002,Ohsawa2004} a Skoda-type $L^2$ division theorem as a corollary of the Ohsawa--Takegoshi $L^2$ extension theorem by reformulating the division problem as an extension problem for the hyperplane bundle on the projectivization of the dual bundles involved in the division problem (see also \cite[Section 3.2.3]{Ohsawa2015}).

The goal of this paper is to show that conversely one can prove a Manivel-type $L^2$ extension theorem as a corollary of a version of the $L^2$ division theorem with bounded generators. %

\begin{definition*}
    Let $X$ be a complex manifold.  If $h$ is a Hermitian metric for a holomorphic vector bundle $V \to X$ with curvature $\Theta_h$, we define the functions $\lambda_h, \Lambda_h : T_X^{1,0} \to \mathbb{R}$ by 
    \[
        \lambda_h(\zeta) := \min_{v \in V_x} \frac{h((\Theta_h)_{\zeta\bar{\zeta}}v,\bar{v})}{h(v,\bar{v})}
        \quad\text{and}\quad
        \Lambda_h(\zeta) := \max_{v \in V_x} \frac{h((\Theta_h)_{\zeta\bar{\zeta}}v,\bar{v})}{h(v,\bar{v})}
    \]
    for all $\zeta \in T_{X,x}^{1,0}$.
\end{definition*}

To simplify notation, denote by $r_V$ the rank of the vector bundle $V \to X$.

\begin{maintheorem}\label{prop:div2ext}
    Let $X$ be an essentially Stein manifold and let $Z = (T=0)$ for some holomorphic section $T$ of a holomorphic vector bundle $E_Z \to X$ of rank $k = \rk E_Z$.  Assume that $E_Z \to X$ carries a Hermitian metric $h_Z$ such that $h_Z(T,\bar{T}) \leq 1$ and that $T$ is generically transverse to the zero section.

    Let $L \to X$ be a holomorphic line bundle with (singular) Hermitian metric $\e^{-\varphi}$ such that 
    \[
        \begin{cases}
            \I\ddbar\varphi \geq (k+2 + \alpha k)\Lambda_{h_Z} + 2\I\ddbar\log\det h_Z \\
            \I\ddbar\varphi \geq 2\I\ddbar\log\det h_Z
        \end{cases}
    \]
    for some $\alpha>1$.  Then for any holomorphic section $f$ of $L|_Z \otimes K_Z \to Z$ such that 
    \[
        \norm{f}_Z^2 := \int_Z |f|^2 \e^{-\varphi} < +\infty
    \]
    there is a holomorphic section $F$ of $L \otimes \det E_Z \otimes K_X \to X$ such that $F|_Z = f \wedge \det(\dif T)$ and 
    \[
        \norm{F}_X^2 := \int_X |F|^2 \e^{-\varphi-\lambda} \leq (k+1) \sigma_k \frac{\alpha}{\alpha-1} \norm{f}_Z^2,
    \]
    where $\sigma_k := \pi^k / k!$ is the volume of the unit ball in $\mathbb{R}^{2k}$.
\end{maintheorem}

A special case is when $Z$ is of codimension $k=1$: if we assume that $T$ is a holomorphic section of some holomorphic line bundle $L_Z \to X$ with metric $\e^{-\eta}$ such that $|T|^2 \e^{-\eta} \leq 1$ and $\dif T|_Z \not\equiv 0$, and we take $\alpha = 1 + \delta$ ($\delta>0$), then the curvature conditions become 
\[
    \begin{cases}
        \I\ddbar\varphi \geq (2 + \delta)\I\ddbar\eta \\
        \I\ddbar\varphi \geq -2\I\ddbar\eta
    \end{cases}
\]
and the constant is $(1+1) \sigma_1 (1 + \delta^{-1}) = 2\pi (1 + \delta^{-1})$.

Even though \autoref{prop:div2ext} is not optimal in both the curvature conditions and the constants involved (see \cite{Manivel1993,Demailly2000,GuanZhou2015}), its proof reduces to a quite straightforward computation once one has the right $L^2$ division theorem.  The main idea is to start with any holomorphic extension $g$ with finite $L^2$ norm, and then ``correct'' it to an extension with smaller norm by using a division theorem (\autoref{prop:bdd-division} below) to remove the holomorphic component of $g$ that is orthogonal to the kernel of the restriction map.

The space of all extensions of a given section on $Z$ is an affine space modeled on the space of sections vanishing on $Z$, i.e.~the ideal generated by the sections that cut out $Z$.  Hence the extension problem becomes a minimization problem in that affine space, meaning an ``(affine) $L^2$ Nullstellensatz'' question, i.e.~$L^2$ division.  Concretely, suppose $X$ is a domain in $\mathbb{C}^n$ and $Z$ is cut out as the zero set of some a holomorphic function $T$.  Given an extension $g$, any other extension has the form $F = g - \beta T$, so that $g = T \beta + 1 F$, a division problem with generators $(T,1)$.  If one rescales $T$ by a very large factor, then one ``distributes'' the $L^2$ norm of the solution $(\beta,F)$ more in the norm of $\beta$ and less in the norm of $F$, minimizing the norm of $F$ in the limit.  From the point of view of the norms, this means that the part of $g$ not along $Z$ contributes less and less when $T$ is rescaled, and in the limit the only remaining contribution is the integral along $Z$, i.e.~the norm of the section to extend.  The precise argument is explained in \autoref{sec:extension}.

It turns out that in order to prove \autoref{prop:div2ext} one needs to be able to take the parameter $\alpha>1$ of the usual Skoda-type theorems to be 1 in the datum norm (see \autoref{rmk:alpha=1}).  This is generally impossible, but if the generators are assumed to be uniformly bounded, then one can prove the following theorem on generically surjective holomorphic morphisms of holomorphic vector bundles.

\begin{maintheorem}\label{prop:bdd-division}
    Let $X$ be an essentially Stein manifold, and let $V,W \to X$ be holomorphic vector bundles with Hermitian metrics $h_V, h_W$, respectively.  Let $V \overset{\gamma}{\longto} W$ be a generically surjective holomorphic morphism of vector bundles such that $(h_V^* \otimes h_W)(\gamma, \bar{\gamma}) \leq 1$.  Fix $\alpha > 1$ and assume that 
    \[
        \begin{cases}
            (r_W + 1) \lambda_{h_W} + (r_V - 1)\Tr\Theta_{h_W} - r_W \Tr\Theta_{h_V} \geq (r_V r_W+1 + \alpha(r_V r_W-1)) (\Lambda_{h_W} - \lambda_{h_V}) \\
            (r_W + 1) \lambda_{h_W} + (r_V - 1)\Tr\Theta_{h_W} - r_W \Tr\Theta_{h_V} \geq (r_V r_W+1) (\Lambda_{h_W} - \lambda_{h_V})
        \end{cases}.
    \]
    Then for any holomorphic section $g$ of $W \otimes K_X \to X$ such that
    \[
        \norm{g}_W^2 := \int_X \frac{h_W(g,\bar{g})}{(h_V^* \otimes h_W)(\gamma,\bar{\gamma})^{r_V r_W}} < +\infty
    \]
    there is a holomorphic section $f$ of $V \otimes K_X \to X$ such that $g = \gamma f$ and
    \[
        \norm{f}_V^2 := \int_X h_V(f,\bar{f}) \leq r_V \frac{\alpha}{\alpha-1} \norm{g}_W^2.
    \]
\end{maintheorem}

Even though the norm of the solution is weaker than the one in the usual Skoda-type theorems, when one has $r_V r_W \leq \dim X$ (such as in the proof of \autoref{prop:div2ext}) \autoref{prop:bdd-division} is, to the author's knowledge, the only Skoda-type theorem where the weight in the norm of the datum is exactly given by the $(-r_V r_W)$-power of the norm of the morphism $\gamma$.  In fact, even quite general statements such as the ones in \cite{Varolin2008} seem only capable of achieving weights growing faster than $(h_V^* \otimes h_W)(\gamma,\bar{\gamma})^{r_V r_W}$ along the zero set of the generators.

In the special case $V = E^{\oplus r}$ and $W = G$ for $E,G \to X$ holomorphic line bundles with metrics $\e^{-\varphi}, \e^{-\psi}$, respectively, one has
\[
    \begin{split}
        \lambda_{h_W} &= \Lambda_{h_W} = \Tr\Theta_{h_W} = \I\ddbar\psi \\
        \lambda_{h_V} &= \frac{1}{r} \Tr\Theta_{h_V} = \I\ddbar\varphi
    \end{split}
\]
Hence we have the following line bundle version of \autoref{prop:bdd-division}.
\begin{maincorollary}\label{prop:bdd-division-lb}
    Let $X$ be a Stein manifold, and let $E,G \to X$ be holomorphic line bundles with (singular) Hermitian metrics $\e^{-\varphi}, \e^{-\psi}$, respectively.  Fix $\gamma_1,\dots,\gamma_r$ holomorphic sections of $E^* \otimes G \to X$ and $\alpha>1$.  Assume that $|\gamma|^2 \e^{-\psi+\varphi} \leq 1$, and that $\I\ddbar\varphi \geq 0$ and
    \[
        \I\ddbar\varphi \geq \frac{\alpha(r-1)}{\alpha(r-1)+1} \I\ddbar\psi.
    \]
    Then for any holomorphic section $g$ of $G \otimes K_X \to X$ such that 
    \[
        \norm{g}_G^2 = \int_X \frac{|g|^2 \e^{-\psi}}{(|\gamma|^2 \e^{-\psi+\varphi})^r}
    \]
    there are holomorphic sections $f_1,\dots,f_r$ of $E \otimes K_X \to X$ such that
    \[
        g = \gamma_1 f_1 + \dots + \gamma_r f_r
    \]
    and 
    \[
        \norm{f}_F^2 = \int_X |f|^2 \e^{-\varphi} \leq r \frac{\alpha}{\alpha-1} \norm{g}_G^2.
    \]
\end{maincorollary}

The proof of \autoref{prop:bdd-division}, in \autoref{sec:division}, is fundamentally the same as \cite{Albesiano2025} and is based on a degeneration argument originally due to B.~Berndtsson and L.~Lempert \cite{BerndtssonLempert2016,Lempert2017}.  Because \autoref{prop:bdd-division} resembles the main result of \cite{Albesiano2025}, and because the latter is not optimal (see also the discussion in \cite{Albesiano2024}), \autoref{prop:bdd-division} is likely also not optimal, both in the curvature conditions and in the constants involved.  %

As an aside, \autoref{prop:bdd-division-lb} can also be used to give a very quick proof of the following Briançon--Skoda-type statement (see \cite[Theorem 11.17]{Demailly2001}).

\begin{maincorollary}\label{prop:BrianconSkoda}
    Fix $h_1, \dots, h_r \in \mathcal{O}_{\mathbb{C}^n,0}$.  Then
    \[
        \mathcal{I}\left( (|h_1|^2 + \dots + |h_r|^2)^{-n} \right) \subset \langle h_1, \dots, h_r \rangle,
    \]
    i.e.~the multiplier ideal sheaf associated to the plurisubharmonic weight $n \log \sum_{j=1}^r|h_j|^2$ is contained in the ideal generated by $h_1, \dots, h_r$.
\end{maincorollary}

The proof of \autoref{prop:BrianconSkoda} is in \autoref{sec:BrianconSkoda}.

\subsection*{Acknowledgements}
I am thankful to Dror Varolin for making me aware of Ohsawa's argument for $L^2$ division via $L^2$ extension, the starting point of what led to the present work, and for many stimulating discussions.  I am also grateful to Bo Berndtsson for precious feedback on an early draft of this paper.

\section{From division to extension}\label{sec:extension}

We illustrate the idea in a simplified setting first.  Take $X = \mathbb{D}^n$, a product of complex disks of radius 1, and let $Z = ( z_1 = 0 )$. Clearly $|z_1|^2 \leq 1$ on $\mathbb{D}^n$.  Let $f$ be a holomorphic function on $Z = \mathbb{D}^{n-1}$ such that 
\[
    \norm{f}_Z^2 := \int_Z |f(z')|^2 \dV(z') < +\infty,
\]
and let $g$ be a holomorphic function on $\mathbb{D}^n$ that is smooth up to the boundary of $\mathbb{D}^n$ and such that $g|_Z = f$ and $\int_{\mathbb{D}^n} |g|^2 \dV < +\infty$.

We apply \autoref{prop:bdd-division-lb} with $E,G$ the trivial line bundles with trivial metrics, and with generators $\gamma_1 := \sqrt{\frac{t}{1+t}} z_1$ and $\gamma_2 := \frac{1}{\sqrt{1+t}}$, for $t \geq 0$ (eventually $t \to +\infty$).  First we check that
\[
    \frac{1}{1+t} \leq |\gamma|^2 = \frac{t}{1+t} |z_1|^2 + \frac{1}{1+t} \leq 1,
\]
so that
\[
    \int_{\mathbb{D}^n} \frac{|g|^2}{|\gamma|^4} \dV \leq (1+t)^2 \int_{\mathbb{D}^n} |g|^2 \dV < +\infty.
\]
(here $r=2$).  Notice also that the curvature conditions are trivially satisfied for all $\alpha > 1$.

\autoref{prop:bdd-division-lb} gives then holomorphic functions $F^{(t)}_1, F^{(t)}_2$ on $\mathbb{D}^n$ such that 
\[
    \sqrt{\frac{t}{1+t}} z_1 F^{(t)}_1 + \frac{1}{\sqrt{1+t}} F^{(t)}_2 = g
\]
and 
\[
    \int_{\mathbb{D}^n} \left( |F^{(t)}_1|^2 + |F^{(t)}_2|^2 \right) \dV \leq \frac{2\alpha}{\alpha-1} \int_{\mathbb{D}^n} \frac{|g|^2}{|\gamma|^4} \dV.
\]

Set now $F^{(t)} := \frac{1}{\sqrt{1+t}} F^{(t)}_2$, so that $F^{(t)}$ is a new extension of $f$.  We have
\[
    \begin{split}
        \norm{F^{(t)}}_{\mathbb{D}^n}^2 &\leq \frac{2\alpha}{\alpha-1} \frac{1}{(1+t)} \int_{\mathbb{D}^n} \frac{|g|^2}{\left( \frac{t}{1+t} |z_1|^2 + \frac{1}{1+t} \right)^2} \dV \\
        &= \frac{2\alpha}{\alpha-1} (1+t) \int_{\mathbb{D}^n} \frac{|g|^2}{\left(t |z_1|^2 + 1 \right)^2} \dV
    \end{split}
\]

Fix $A>1$ and let $\varepsilon>0$ be such that $|g(z_1,z')|^2 \leq A |f(z')|^2$ for all $|z_1| < \varepsilon$ and $z' \in \mathbb{D}^{n-1}$.  Then 
\[
    \int_{\mathbb{D}^n} \frac{|g|^2}{\left(t |z_1|^2 + 1 \right)^2} \dV \leq A \int_{|z_1| < \varepsilon} \frac{|f|^2}{\left(t |z_1|^2 + 1 \right)^2} \dV + \int_{|z_1| \geq \varepsilon} \frac{|g|^2}{(t \varepsilon + 1)^2} \dV.
\]
The second term decays like $t^{-2}$ as $t \to +\infty$, so that 
\[
    \begin{split}
        \lim_{t \to +\infty} \norm{F^{(t)}}_{\mathbb{D}^n}^2 &\leq \frac{2\alpha A}{\alpha-1} \lim_{t \to +\infty} (1+t) \int_{|z_1| < \varepsilon} \frac{|f|^2}{\left(t |z_1|^2 + 1 \right)^2} \dV \\
        &= \frac{2\alpha A}{\alpha-1} \norm{f}_Z^2 \lim_{t \to +\infty} (1+t) \int_{|z_1| < \varepsilon} \frac{\dif z_1 \wedge \dif \bar{z}_1}{\left(t |z_1|^2 + 1 \right)^2} \\
        &= 2\pi \frac{\alpha A}{\alpha-1} \norm{f}_Z^2 \lim_{t \to +\infty} (1+t) \int_0^\varepsilon \frac{2\rho \dif \rho}{(t \rho^2 + 1)^2} \\
        &= 2\pi \frac{\alpha A}{\alpha-1} \norm{f}_Z^2 \lim_{t \to +\infty} \left(\frac{1}{t}+1\right) \left[ 1 - \frac{1}{t \varepsilon^2 + 1} \right] = 2\pi \frac{\alpha A}{\alpha-1} \norm{f}_Z^2
    \end{split}
\]
for all $A > 1$.  We then get a sequence of extensions $F^{(t)}$ (also depending on $\alpha>1$) such that
\[
    \lim_{\substack{t \to +\infty \\ \alpha \to +\infty}} \norm{F^{(t)}}_{\mathbb{D}^n} \leq 2\pi \norm{f}_Z^2.
\]
By the sub-mean value property and Montel's theorem we then obtain a holomorphic function $F$ extending $f$ and such that $\norm{F}_{\mathbb{D}^n}^2 \leq 2\pi \norm{f}_Z^2$.

\begin{remark}\label{rmk:alpha=1}
    If we use the standard theorem of Skoda we achieve the following estimate:
    \[
        \begin{split}
            \lim_{t \to +\infty} \norm{F^{(t)}}_{\mathbb{D}^n}^2 &\leq \frac{2\alpha A}{\alpha-1} \lim_{t \to +\infty} (1+t)^\alpha \int_{|z_1| < \varepsilon} \frac{|f|^2}{\left(t |z_1|^2 + 1 \right)^{\alpha+1}} \dV \\
            &= 2\pi \frac{\alpha A}{\alpha-1} \norm{f}_Z^2 \lim_{t \to +\infty} (1+t)^\alpha \int_0^\varepsilon \frac{2\rho \dif \rho}{(t \rho^2 + 1)^{\alpha+1}} \\
            &= 2\pi \frac{\alpha A}{\alpha-1} \norm{f}_Z^2 \lim_{t \to +\infty} \frac{(1+t)^\alpha}{\alpha t} \left[ 1 - \frac{1}{\left(t \varepsilon^2 + 1\right)^\alpha} \right] = +\infty
        \end{split}
    \]
    A similar computation shows that the bound blows up when $t \to +\infty$ if instead we use the $L^2$ division theorems with bounded generators of \cite{Varolin2008}, because of the extra logarithmic term present in the norms of Theorems 2.2 and 2.3 there.

    Hence, it really seems necessary to be able to take the power on the weight at the denominator of the datum norm in the division problem to be exactly $\codim Z + 1$ in the way achieved by \autoref{prop:bdd-division}.
\end{remark}

\begin{proof}[Proof of \autoref{prop:div2ext}]
As in \cite{BerndtssonLempert2016} and \cite{Albesiano2024}, because of the universality of the bounds, by shrinking $X$ and removing the locus where $Z$ is singular we can assume that $X$ is Stein, that $\det \dif T|_Z \neq 0$ everywhere, and that $f$ is smooth up to the boundary of $Z$. We can also assume that the metrics involved are smooth, and that there is an extension $g$ smooth up to the boundary of $X$ and with finite $L^2$-norm, namely: there is a holomorphic section $g$ of $L \otimes \det E_Z \otimes K_X \to X$, smooth in a neighborhood of $\bar{X}$, such that $g|_Z = f \wedge \det(\dif T)$ and $\norm{g}_X^2 < +\infty$.

The idea is then to use \autoref{prop:bdd-division} to ``correct'' $g$ to an extension $F$ with smaller norm.  Specifically, we will apply the division theorem with $V := (E_Z^* \oplus \mathcal{O}_X) \otimes L \otimes \det E_Z$ and $W := L \otimes \det E_Z$.  Denote by $\e^{-\eta} := \det h_Z$ the metric on $\det E_Z$ induced by $h_Z$.  We endow $V$ with the Hermitian metric
\[
    h_V = \begin{pmatrix}
        h_Z^* \e^{-\varphi-\eta} & 0 \\ 0 & \e^{-\varphi-\eta}
    \end{pmatrix}
\]
and $W$ with the metric $h_W = \e^{-\varphi-\eta}$.

In this setting we have
\[
    \begin{split}
        r_V &= k+1, \quad r_W = 1, \\
        \lambda_{h_W} &= \Lambda_{h_W} = \Tr\Theta_{h_W} = \I\ddbar(\varphi+\eta), \\
        \lambda_{h_V} &= \min(-\Lambda_{h_Z}, 0) + \I\ddbar(\varphi+\eta), \\
        \Tr\Theta_{h_V} &= -\Tr\Theta_{h_Z} + (k+1) \I\ddbar(\varphi+\eta) = (k+1)\I\ddbar\varphi + k\I\ddbar\eta.
    \end{split}
\]
Then 
\[
    (r_W + 1) \lambda_{h_W} + (r_V - 1)\Tr\Theta_{h_W} - r_W \Tr\Theta_{h_V} = \I\ddbar(\varphi+2\eta)
\]
and
\[
    \Lambda_{h_W} - \lambda_{h_V} = - \min(-\Lambda_{h_Z},0) = \max(\Lambda_{h_Z},0),
\]
so that in order to apply \autoref{prop:bdd-division} we need 
\[
    \begin{split}
    \I\ddbar\varphi &\geq (k+2 + \alpha k)\Lambda_{h_Z} - 2\I\ddbar\eta \\
        &= (k+2 + \alpha k)\Lambda_{h_Z} + 2\I\ddbar\log\det h_Z
    \end{split}
\]
and
\[
    \I\ddbar\varphi \geq -2\I\ddbar\eta = 2\I\ddbar\log\det h_Z.
\]

Fix now $t \geq 0$.  We see $T$ as a holomorphic morphism of vector bundles 
\[
    E_Z^* \otimes L \otimes \det E_Z \, \overset{T}{\longto} \, L \otimes \det E_Z
\]
and take $V \overset{\gamma}{\longto} W$ to be given by
\[
    \gamma(v_1,v_2) := \sqrt{\frac{t}{1+t}} T v_1 + \frac{1}{\sqrt{1+t}} v_2.
\]

Note that by construction
\[
    \frac{1}{1+t} \leq (h_V^* \otimes h_W)(\gamma,\bar{\gamma}) = \frac{t}{1+t} h_Z(T,\bar{T}) + \frac{1}{1+t} \leq 1,
\]
so that 
\[
    \norm{g}_W^2 = \int_X \frac{|g|^2 \e^{-\varphi-\eta}}{(h_V^* \otimes h_W)(\gamma,\bar{\gamma})^{k+1}} \leq (1+t)^{k+1} \int_X |g|^2 \e^{-\varphi-\eta} < +\infty.
\]

We can then apply \autoref{prop:bdd-division} to obtain $F_1^{(t)} \in H^0(X, E_Z^* \otimes L \otimes \det E_Z \otimes K_X)$ and $F_2^{(t)} \in H^0(X, L \otimes \det E_Z \otimes K_X)$ such that 
\[
    \sqrt{\frac{t}{1+t}} T F_1^{(t)} + \frac{1}{\sqrt{1+t}} F_2^{(t)} = g
\]
and
\[
    \norm{(F_1^{(t)},F_2^{(t)})}_V^2 = \int_X \left(h_Z^*(F_1^{(t)},\overline{F}_1^{(t)}) + |F_2^{(t)}|^2\right) \e^{-\varphi-\eta} \leq \frac{\alpha(k+1)}{\alpha-1} \norm{g}_W^2.
\]

The holomorphic section $F^{(t)} := (1+t)^{-1/2} F_2^{(t)}$ of $L \otimes \det E_Z \otimes K_X \to X$ provides now a new extension of $f$:
\[
    F^{(t)}|_Z = \frac{F_2^{(t)}|_Z}{\sqrt{1+t}} = g|_Z - \sqrt{\frac{t}{1+t}} (T F_1^{(t)})|_Z = g|_Z = f \wedge \det(\dif T).
\]
Moreover
\[
    \begin{split}
        \norm{F^{(t)}}_X^2 &= \int_X |F^{(t)}|^2 \e^{-\varphi-\eta} = \frac{1}{1+t} \int_X |F_2^{(t)}|^2 \e^{-\varphi-\eta} \leq \frac{1}{1+t} \norm{(F_1^{(t)},F_2^{(t)})}_W^2 \\
        &\leq \frac{1}{1+t} \frac{\alpha(k+1)}{\alpha-1} \int_X \frac{|g|^2 \e^{-\varphi-\eta}}{\left( \frac{t}{1+t}h_Z(T,\bar{T}) + \frac{1}{1+t} \right)^{k+1}} \\
        &= \frac{\alpha(k+1)}{\alpha-1} (1+t)^k \int_X \frac{|g|^2 \e^{-\varphi-\eta}}{\left( t h_Z(T,\bar{T}) + 1 \right)^{k+1}}.
    \end{split}
\]

Since the extension $g$ and the metrics involved are smooth up to the boundary, for every $A>1$ there is $\varepsilon>0$ such that $\mset{h_Z(T,\bar{T}) < \varepsilon^2}$ is a tube around $Z$ in which 
\[
    |g|^2 \e^{-\varphi-\eta} \leq A |\tilde{f}|^2 \e^{-\varphi} \wedge |\!\det \dif T|^2 \e^{-\eta},
\]
where $\tilde{f}$ is the pullback of $f$ via the projection of $\mset{h_Z(T,\bar{T}) < \varepsilon^2}$ to $Z$.  Then 
\[
    \int_X \frac{|g|^2 \e^{-\varphi-\eta}}{\left( t h_Z(T,\bar{T}) + 1 \right)^{k+1}} \leq A \int_{h_Z(T,\bar{T}) < \varepsilon^2} \frac{|\tilde{f}|^2 \e^{-\varphi} \wedge |\!\det\dif T|^2 \e^{-\eta}}{(t h_Z(T,\bar{T})+1)^{k+1}} + \int_{h_Z(T,\bar{T}) \geq \varepsilon^2} \frac{|g|^2 \e^{-\varphi-\eta}}{(t \varepsilon^2+1)^{k+1}}.
\]

Note now that
\[
    \lim_{t \to +\infty} (1+t)^k \int_{h_Z(T,\bar{T}) \geq \varepsilon^2} \frac{|g|^2 \e^{-\varphi-\eta}}{(t \varepsilon^2+1)^{k+1}} \leq \lim_{t \to +\infty} \frac{(1+t)^k}{(1+\varepsilon^2 t)^{k+1}} \norm{g}_Z^2 = 0
\]
for all $\varepsilon>0$.  We compute the other term as follows:
\[
    \begin{split}
        &\int_{h_Z(T,\bar{T}) < \varepsilon^2} \frac{|\tilde{f}|^2 \e^{-\varphi} \wedge |\!\det\dif T|^2 \e^{-\eta}}{(t h_Z(T,\bar{T}) + 1)^{k+1}} \\
        &\quad= \left(\,\norm{f}_Z^2 + O(\varepsilon)\right) \frac{2 \pi^k}{(k-1)!} \int_0^\varepsilon \frac{\rho^{2k-1}}{(t \rho^2 + 1)^{k+1}} \dif \rho \\
        &\quad= \frac{\pi^k}{k!} \left(\,\norm{f}_Z^2 + O(\varepsilon)\right) \frac{\varepsilon^{2k}}{(t \varepsilon^2 + 1)^k},
    \end{split}
\]
so that 
\[
    \begin{split}
        &\lim_{t \to +\infty} A(1+t)^k \int_{h_Z(T,\bar{T}) < \varepsilon^2} \frac{|\tilde{f}|^2 \e^{-\varphi} \wedge |\!\det\dif T|^2 \e^{-\eta}}{(t h_Z(T,\bar{T}) + 1)^{k+1}} \\
        &\quad= \frac{\pi^k}{k!} A\left(\,\norm{f}_Z^2 + O(\varepsilon)\right) \lim_{t \to +\infty} \frac{(1+t)^k \varepsilon^{2k}}{(t\varepsilon^2 + 1)^k} = A \sigma_k \left(\,\norm{f}_Z^2 + O(\varepsilon)\right)
    \end{split}
\]
for all $A>1$, where $\sigma_k := \pi^k/k!$ is the volume of the unit ball in $\mathbb{R}^{2k}$.  Altogether we obtain
\[
    \lim_{t \to +\infty} (1+t)^k \int_X \frac{|g|^2 \e^{-\varphi-\eta}}{\left( t h_Z(T,\bar{T}) + 1 \right)^{k+1}} \leq A \sigma_k \left(\,\norm{f}_Z^2 + O(\varepsilon)\right)
\]
for all $A > 1$.  By taking the limit $A \to 1$ (and thus $\varepsilon \to 0$) we conclude 
\[
    \lim_{t \to +\infty} \norm{F^{(t)}}_X^2 \leq (k+1) \sigma_k \frac{\alpha}{\alpha-1} \norm{f}_Z^2.
\]

Overall we then have a sequence of holomorphic sections $F^{(t)}$ of $L \otimes L_Z \otimes K_X \to X$ extending $f$ (i.e.~such that $F^{(t)}|_Z = f \wedge \dif T$) and with an upper bound on the $L^2$-norms converging to the $L^2$-norm of the datum $f$, up to a constant.  By the sub-mean value property, uniform estimates with plurisubharmonic weights imply locally uniform sup-norm estimates.  Then by Montel's theorem $F^{(t)}$ converges up to subsequences to a holomorphic section $F$ of $L \otimes \det E_Z \otimes K_X \to X$.  Clearly $F|_Z = f \wedge \det(\dif T)$ and 
\[
    \norm{F}_X^2 \leq (k+1) \sigma_k \frac{\alpha}{\alpha-1} \norm{f}_Z^2,
\]
proving \autoref{prop:div2ext}.
\end{proof}

\section{A division theorem with bounded generators (proof of \autoref{prop:bdd-division})}\label{sec:division}
The proof of \autoref{prop:bdd-division} follows \cite{Albesiano2025} almost step by step.  First, because the constant $r_V \frac{\alpha}{\alpha-1}$ does not depend on $X$, we can work under the assumption that $\gamma$ is actually surjective by removing the divisor of $\gamma^{\wedge r_G}$, and then extending across the divisor by Riemann's theorem on removable singularities.  We can also assume that $X$ is Stein.

The whole point of the argument is then to find a positively curved family of norms that interpolate between the norm $\norm{\cdot}_W$ of the datum and the norm $\norm{\cdot}_V$ of the solution, the rest of the argument following by the construction of Berndtsson and Lempert \cite{BerndtssonLempert2016} based on Berndtsson's theorem on direct image bundles \cite{Berndtsson2009} (here in its vector bundle version \cite[Theorem 2]{Varolin2025}).  Denote by $\pr_X: P(V^* \otimes W) \to X$ the bundle projection, set $r := \rk(V^* \otimes W)= r_V r_W$, and set $h := h_V^* \otimes h_W$.  We work with the vector bundle 
\[
    \pr_X^* (W \otimes K_X) \otimes \mathcal{O}(1) \longto P(V^* \otimes W),
\]
endowed with the following family of metrics, parametrized by 
\[
    \tau \in \mathbb{L} := \{ \Re z \leq 0 \}.
\]
Recall first that a section $s$ of $V \to X$ lifts to the section $\hat{s}(v) := vs$ of $\pr_X^* (W \otimes K_X) \otimes \mathcal{O}(1) \to P(V^* \otimes W)$.

For a section $\sigma$ of $\pr_X^* (W \otimes K_X) \otimes \mathcal{O}(1) \to P(V^* \otimes W)$ set
\[
    \mathfrak{h}_\tau(\sigma,\bar{\sigma}) := \frac{r_V}{\Vol_{\Re\tau}} \frac{h_W(\sigma_v, \bar{\sigma}_v)}{h(v,\bar{v})} \left( h(v,\bar{v}) \e^{-\chi_\tau} \right)^{\alpha(r-1)},
\]
where $\Vol_t = \frac{\pi^{r-1}}{(r-1)!} \e^{(r-1)t}$ is the volume of the ball of radius $\e^{t/2}$ in $\mathbb{R}^{2r-2}$ and
\[
    \chi_\tau := \max\left( \log\left[h(\gamma,\bar{\gamma}) h(v,\bar{v}) - |h(\gamma,\bar{v})|^2\right] - \Re\tau, \, \log h(v,\bar{v}) \right).
\]
Note that, compared to the weight in \cite{Albesiano2025}, the right-hand side of the maximum does not contain $h(\gamma,\bar{\gamma})$.

Since the general outline of the proof is the same as \cite{Albesiano2025}, we only highlight the differences, namely the endpoints of the family and the new curvature conditions.

\subsection{The family of norms}
Let $\omega := \I\ddbar\log h(v,\bar{v})$ and set 
\[
    \norm{\sigma}_\tau^2 := \frac{r_V}{\Vol_{\Re\tau}} \int_{[v] \in P(V^* \otimes W)} \frac{h_W(\sigma_v, \bar{\sigma}_v)}{h(v,\bar{v})} \left( h(v,\bar{v}) \e^{-\chi_\tau} \right)^{\alpha(r-1)} \wedge (\iota_x^*\omega)^{\wedge (r - 1)}
\]
(here and in the following $x = \pr_X(v)$).  Since $\chi_\tau$, $\mathfrak{h}_\tau$, and $\norm{\cdot}_\tau$ clearly only depend on $t = \Re\tau$, we will write $\chi_t$, $\mathfrak{h}_t$, and $\norm{\cdot}_t$, respectively, with the understanding that $t \in (-\infty,0]$.

Because
\[
    h(\gamma,\bar\gamma) h(v,\bar{v}) - |h(\gamma,\bar{v})|^2 \leq h(\gamma,\bar\gamma)h(v,\bar{v}) \leq h(v,\bar{v}),
\]
when $t=0$ we have $\chi_0 = \log h(v,\bar{v})$.  Hence
\[
    \norm{\hat{s}}_0^2 = \frac{r_V}{\Vol_0} \int_{v \in P(V^* \otimes W)} \frac{h_W(vs,\overline{vs})}{h(v,\bar{v})} \wedge (\iota_x^*\omega)^{\wedge (r - 1)}.
\]
The same computation of \cite[Section 2]{Albesiano2025} gives then $\norm{\hat{s}}_0^2 = \norm{s}_V^2$.

Next we look at $t \to -\infty$.  Fix $x \in X$ and $t < 0$, and define $A_{t,x}$ to be the subset of the fiber $P(V^* \otimes W)_x$ where the maximum in $\chi_t$ is attained by the right-hand side, i.e.
\[
    A_{t,x} := \mset{v \in P(V^* \otimes W)_x}[1 - \frac{|h(\gamma,\bar{v})|^2}{h(\gamma,\bar\gamma)h(v,\bar{v})} \leq \frac{\e^t}{h(\gamma,\bar\gamma)}].
\]
Then $A_{t,x}$ is a ball of real dimension $2r-2$ centered at $[h(x)]$ and with radius asymptotic to $\sqrt{\frac{\e^t}{h(\gamma,\bar{\gamma})}}$ when $t \to -\infty$.

We then write $\norm{\hat{s}}_t^2 = \one_t + \two_t$ with 
\[
    \one_t := r_V \fint_{A_t} \frac{h_W(vs, \overline{vs})}{h(v,\bar{v})h(\gamma,\bar{\gamma})^{r-1}} \wedge (\iota_x^*\omega)^{\wedge (r - 1)}
\]
and
\[
    \two_t := \frac{r_V}{\Vol_t} \int\limits_{x \in X} \int\limits_{[v] \in P(V^* \otimes W)_x \setminus A_{t,x}} \frac{h_W(vs, \overline{vs}) \e^{\alpha(r-1)t}}{h(v,\bar{v}) h(\gamma,\bar\gamma)^{\alpha(r-1)}} \wedge \frac{(\iota_x^*\omega)^{\wedge (r - 1)}}{\left( 1 - \frac{|h(\gamma,\bar{v})|^2}{h(\gamma,\bar{\gamma})h(v,\bar{v})}\right)^{\alpha(r-1)}}.
\]

Integrating along the fibers $P(V^* \otimes W)_x$ first, one obtains
\begin{equation}\label{eq:It}
    \lim_{t \to -\infty} \one_t = r_V \int_X \frac{h_W(\gamma s, \overline{\gamma s})}{h(\gamma,\bar{\gamma})^r}.
\end{equation}
The second term can be rewritten as 
\[
    \two_t = \e^{-(r-1)t} \int_t^0 \e^{-\alpha(r-1)(\tilde{t}-t)} \dif\nu_s(\tilde{t})
\]
with $\nu_s(t) = \e^{(r-1)t} \one_t$.  By the calculus lemma of \cite[Lemma 3.1]{Albesiano2024} we then have 
\[
    \lim_{t \to -\infty} \two_t = \frac{r_V}{\alpha-1} \int_X \frac{h_W(\gamma s, \overline{\gamma s})}{h(\gamma,\bar{\gamma})^r},
\]
so that overall 
\begin{equation}\label{eq:limit}
    \lim_{t \to -\infty} \norm{\hat{s}}_t^2 = r_V \frac{\alpha}{\alpha-1} \int_X \frac{h_W(\gamma s, \overline{\gamma s})}{h(\gamma,\bar{\gamma})^r} = r_V \frac{\alpha}{\alpha-1} \norm{\gamma s}_W^2,
\end{equation}
retrieving a multiple of the norm squared of the image of $s$ under $\gamma$.

\subsection{Curvature}
The computations on curvature are essentially the same as Section 3 of \cite{Albesiano2025}, except that $\chi^{(2)}$ is now replaced by $\chi^{(2)} = \log h(v,\bar{v})$.  Condition (3.3) of \cite{Albesiano2025} is therefore the same, while condition (3.4) becomes
\[
    \begin{split}
        (r_W &+ 1)\Theta_{h_W} + (r+1) \I\ddbar\log h(v,\bar{v}) \\
        &+ (r_V - 1) \Theta_{\det h_W} - r_W \Theta_{\det h_V} \geqGrif 0.
    \end{split}
\]
Consequently, (3.5) remains the same and (3.6) becomes
\[
    (r_W + 1)\frac{h_W(\Theta_{h_W}u, \bar{u})}{h_W(u,\bar{u})} + (r_V - 1) \Theta_{\det h_W} - r_W \Theta_{\det h_V} \geq (r+1) \frac{h(\Theta_h v, \bar{v})}{h(v,\bar{v})},
\]
meaning that it suffices to have
\[
    \begin{split}
        (r_W + 1) &\lambda_{h_W} + (r_V - 1) \Tr\Theta_{h_W} - r_W \Tr\Theta_{h_V} \\
        &\geq \alpha(r-1) \Lambda_{h^{\otimes2}} + (r+1 - \alpha(r-1)) \Lambda_h = (r+1 + \alpha(r-1)) \Lambda_h
    \end{split}
\]
for the left-hand side of max, and
\[
    (r_W + 1) \lambda_{h_W} + (r_V - 1) \Tr\Theta_{h_W} - r_W \Tr\Theta_{h_V} \geq(r+1) \Lambda_h
\]
for the right-hand side.  Since
\[
    \Lambda_h = \Lambda_{h_F^* \otimes h_G} = \Lambda_{h_F^*} + \Lambda_{h_G} = \Lambda_{h_G} - \lambda_{h_F}
\]
(see \cite[Section 1.1]{Albesiano2025}) we conclude that the requirements on curvature are satisfied when
\[
    \begin{cases}
        (r_W + 1) \lambda_{h_W} + (r_V - 1)\Tr\Theta_{h_W} - r_W \Tr\Theta_{h_V} \geq (r_V r_W+1 + \alpha(r_V r_W-1)) (\Lambda_{h_W} - \lambda_{h_V}) \\
        (r_W + 1) \lambda_{h_W} + (r_V - 1)\Tr\Theta_{h_W} - r_W \Tr\Theta_{h_V} \geq (r_V r_W+1) (\Lambda_{h_W} - \lambda_{h_V})
    \end{cases},
\]
as in the hypotheses of \autoref{prop:bdd-division}.

\section{A Briançon--Skoda-type theorem (proof of \autoref{prop:BrianconSkoda})}\label{sec:BrianconSkoda}
Let $\mathcal{I}$ be the ideal generated by the germs $h_1, \dots, h_r \in \mathcal{O}_{\mathbb{C}^n,0}$, and let $\mathcal{J}$ be the multiplier ideal sheaf associated to the weight $n\log(\sum_j|h_j|^2)$, namely
\[
    \mathcal{J} := \mset{g \in \mathcal{O}_{\mathbb{C}^n,0}}[\frac{|g|^2}{(|h_1|^2 + \dots + |h_r|^2)^n} \in L^1_\text{loc}].
\]

We recall the following lemma \cite[Lemma 11.16]{Demailly2001}.
\begin{lemma*}
    If $\mathcal{I} = \langle h_1, \dots, h_r \rangle \subset \mathcal{O}_{\mathbb{C}^n,0}$ and $p > n$, then there are $\tilde{h}_1, \dots, \tilde{h}_n \in \mathcal{I}$ and $C>0$ such that 
    \[
        C^{-1} \sum_{j=1}^r |h_j|^2 \leq \sum_{k=1}^n |\tilde{h}_k|^2 \leq C \sum_{j=1}^r |h_j|^2.
    \]
    In fact, a generic tuple has this property.
\end{lemma*}

Take then $g \in \mathcal{J}$.  We can always assume that $r \geq n$ by setting $h_{r+1} = \dots = h_n = 0$.  Let $\tilde{h}_1, \dots, \tilde{h}_n$ be given by the lemma.  Take also $U$ to be a pseudoconvex neighborhood of $0 \in \mathbb{C}^n$ where $h_1, \dots, h_r$, $\tilde{h}_1, \dots, \tilde{h}_n$, and $g$ are defined, and where $\frac{|g|^2}{(|h_1|^2 + \dots + |h_r|^2)^n}$ is $L^1$.  Apply \autoref{prop:bdd-division-lb} with $E$ and $G$ the trivial bundle with trivial metrics, and with generators $\tilde{h}_1, \dots, \tilde{h}_n$.  We have 
\[
    \norm{g}_W^2 = \int_U \frac{|g|^2}{(|\tilde{h}_1|^2 + \dots + |\tilde{h}_n|^2)^n} \dV \leq C^n \int_U \frac{|g|^2}{(|h_1|^2 + \dots + |h_r|^2)^n} \dV < +\infty,
\]
so we obtain holomorphic functions $f_1, \dots, f_n$ on $U$ such that $g = \tilde{h}_1 f_1 + \dots + \tilde{h}_n f_n$.  Thus $f \in \langle \tilde{h}_1, \dots, \tilde{h}_n \rangle \subset \mathcal{I}$, proving \autoref{prop:BrianconSkoda}.

\end{document}